This article has been accepted for inclusion in a future issue of this journal. Content is final as presented, with the exception of pagination.

# Minimization of Fraction Function Penalty in Compressed Sensing

Haiyang Li, Qian Zhang, Angang Cui, and Jigen Peng

***Abstract*— In this paper, we study the minimization problem of a non-convex sparsity-promoting penalty function, i.e., fraction function, in compressed sensing. First, we discuss the equivalence of $\ell_0$ minimization and fraction function minimization. It is proved that the optimal solution to fraction function minimization solves $\ell_0$ minimization and the optimal solution to the regularization problem also solves fraction function minimization if the certain conditions are satisfied, which is similar to the regularization problem in a convex optimization theory. Second, we study the properties of the optimal solution to the regularization problem, including the first-order and second-order optimality conditions and the lower and upper bounds of the absolute value for its nonzero entries. Finally, we derive the closed-form representation of the optimal solution to the regularization problem and propose an iterative $FP$ thresholding algorithm to solve the regularization problem. We also provide a series of experiments to assess the performance of the $FP$ algorithm, and the experimental results show that the $FP$ algorithm performs well in sparse signal recovery with and without measurement noise.**

## I. Introduction

THE goal of compressed sensing (see [2], [11]) is to reconstruct a sparse signal under a few linear measurements that are far less than the dimension of the ambient space of the signal. The following minimization is commonly employed to model this problem:

$$(P_0)\ \min_{x\in\mathbb{R}^n} \|x\|_0 \quad \text{s.t. } Ax = b \tag{1}$$

where $A$ is an $m \times n$ real matrix of full row rank with $m < n$, $b$ is a nonzero real vector of $m$ dimension, and $\|x\|_0$ is the so-called $\ell_0$-norm of real vector $x$, which counts the number of the nonzero entries in $x$ (see [1], [16], [29]). Sparsity problems can be frequently transformed into the following so-called $\ell_0$ regularization problem:

$$\left(P_0^{\lambda}\right) \min_{x\in\mathbb{R}^n} \|Ax - b\|_2^2 + \lambda\|x\|_0 \tag{2}$$

where $\lambda > 0$, called the regularized parameter, represents a tradeoff between error and sparsity. Thiao [30] shows that there exists $\lambda_0 > 0$, such that the minimization problems $(P_0^{\lambda})$ and $(P_0)$ have the same solution set for all $0 < \lambda \leq \lambda_0$. Unfortunately, although the $\ell_0$-norm characterizes the sparsity of the vector $x$, the $\ell_0$ optimization problem is actually NP-hard because of the discrete and discontinuous nature of the $\ell_0$-norm. In general, the relaxation methods replace the $\ell_0$-norm by continuous sparsity promoting penalty functions $P(\cdot)$. The minimization takes the form

$$\min_{x\in\mathbb{R}^n} P(x) \quad \text{s.t. } Ax = b \tag{3}$$

for the constrained problem and

$$\min_{x\in\mathbb{R}^n} \|Ax - b\|_2^2 + \lambda P(x) \tag{4}$$

for the regularization problem. Convex relaxation uniquely selects $P(x)$ as the $\ell_1$-norm. Considerable excellent theoretical work (see [3], [13]–[15], [21]), together with some empirical evidence (see [6]), has shown that, provided some conditions are met, such as assuming the restricted isometric property (RIP) [3], the $\ell_1$-norm minimization can truly make an exact recovery. According to the convex optimization theory, there exists some $\lambda > 0$, such that the solution to the regularization problem (4) also solves the constrained problem (3) when $P(x) = \|x\|_1$. The $\ell_1$ algorithms for solving the regularization problem include $\ell_1$-magic [2], the soft thresholding algorithm (soft algorithm in brief) (see [8], [12]), the Bregman and split Bregman methods (see [20], [35]), and the alternating direction algorithms [34].

There are many choices of $P(x)$ for non-convex relaxation, in which the $\ell_p$-norm [$p \in (0, 1)$] appears to be the most popular choice. Key work by Gribonval and Nielsen [21] on $0 < p < 1$ has resulted in the above-described optimization models gaining in popularity in the literature (see [5], [7], [10], [17], [18], [22], [23], [26]–[28], [31], [33]). In [26], we have demonstrated that in every underdetermined linear system $Ax = b$, there corresponds a constant $p^*(A, b) > 0$, which is called the $NP/CMP$ equivalence constant ($NP/CMP$ equivalence means that the NP-hard optimization problem is equivalent to the continuous minimization problem), such that every solution to the $\ell_p$-norm minimization problem

Manuscript received August 14, 2018; revised December 26, 2018, April 27, 2019, and June 2, 2019; accepted June 2, 2019. This work was supported in part by the NSFC under Grant 11271297 and Grant 11131006 and in part by the National Basic Research Program of China under Grant 2013CB329404. *(Corresponding author: Jigen Peng.)*

H. Li and J. Peng are with the School of Mathematics and Information Science, Guangzhou University, Guangzhou 510006, China (e-mail: fplihaiyang@126.com; jgpeng@mail.xjtu.edu.cn).

Q. Zhang is with the School of Computer Science and Engineering, University of Electronic Science and Technology, Chengdu 611731, China (e-mail: hnnyzhqian@163.com).

A. Cui is with the School of Mathematics and Statistics, Xi'an Jiaotong University, Xi'an 710049, China (e-mail: cuiangang@stu.xjtu.edu.cn).

also solves the $\ell_0$-norm minimization problem whenever $0 < p < p^*(A, b)$. At present, there are two main algorithmic approaches to $\ell_p$-norm minimization for $0 < p < 1$. One is the iteration reweighted least squares minimization algorithm (the IRLS algorithm in brief) [9]. The authors proved that the rate of local convergence of this algorithm was superlinear and the rate was faster for smaller $p$ and increased toward quadratic as $p \to 0$. Moreover, at each iteration, the solution of a least squares problem is required, of which the computational complexity is $\mathcal{O}(mn^2)$. The other approach is an iterative thresholding algorithm when $p = (1/2), (2/3)$ (see [4], [33]). The authors showed that $\ell_{1/2}$ regularization could be quickly solved by the iterative half thresholding algorithm (the half algorithm in brief) and the algorithm was convergent when applied to the $k$-sparsity problem. In addition, at every iteration step of the half algorithm, some products between the matrix and the vector are required; thus, the computational complexity is $\mathcal{O}(mn)$.

Although the computational complexity of the half algorithm is lower than that of IRLS, we do not know whether there is any $\lambda > 0$, such that the optimal solution to the regularization problem (4) also solves the constrained problem (3) when $P(x) = \|x\|_{0.5}^{0.5}$, which is different from the result when $P(x) = \|x\|_1$.

In this paper, inspired by the good performance of the fraction function $p_a(x) = a|x|/1 + a|x|$, called "strictly non-interpolating" in [19], in image restoration, we take

$$P(x) = P_a(x) = \sum_{i=1}^{n} p_a(x_i) \ (x \in \mathbb{R}^n).$$

In fact, the fraction function is widely used in image restoration. Geman and Reynolds [19] showed that the fraction function gave rise to a step-shaped estimate from the ramp-shaped data. Furthermore, Nikolova [24] demonstrated that for almost all data, the strongly homogeneous zones recovered by the fraction function were preserved constant under any small perturbation of the data. We will study the following minimization problems $(FP_a)$ and $(FP_a^\lambda)$ in terms of theory, algorithms, and computation. The constrained fraction function minimization version is

$$(\mathrm{FP}_a) \min_{x \in \mathbb{R}^n} P_a(x) \quad \text{s.t. } Ax = b \tag{5}$$

and the unconstrained fraction function regularization version is

$$(\mathrm{FP}_a^\lambda) \qquad \min_{x \in \mathbb{R}^n} \left\{ \|Ax - b\|_2^2 + \lambda P_a(x) \right\}. \tag{6}$$

This paper is organized as follows. In Section II, we study the elementary properties of the fraction function. In Section III, we focus on proving some theorems, which establish the equivalence of $(FP_a)$ and $(P_0)$. Section IV is devoted to discussing the equivalence of $(FP_a^\lambda)$ and $(FP_a)$ and the properties of the optimal solution to the regularization problem $(FP_a^\lambda)$, including the first-order and second-order optimality conditions and the lower and upper bounds of the absolute value for its nonzero entries. In Section V, we derive the closed-form representation of the optimal solution to the regularization problem $(FP_a^\lambda)$ by using the Cardano formula on roots of cubic polynomials and algebraic identities and propose an iterative $FP$ thresholding algorithm to solve the regularization problem $(FP_a^\lambda)$. In Section VI, we present the experiments with a series of sparse signal recovery applications to demonstrate the robustness and effectiveness of the new algorithms. We conclude this paper in Section VII.

## II. Preliminaries and the Properties of the Fraction Function

We consider the fraction function

$$p_a(t) = \frac{a|t|}{1 + a|t|}$$

where the parameter $a \in (0, +\infty)$. It is easy to verify that $p_a(t)$ is symmetric, $p_a(t) = 0$ if $t = 0$, and $\lim_{a \to +\infty} p_a(t) = 1$ if $t \neq 0$. Clearly, with the adjustment of parameter $a$, $p_a(t)$ can approximate $\ell_0$ well. Moreover, $p_a(t)$ is increasing and concave for $t \in [0, +\infty)$.

It is also easy to check that the triangle inequality holds for the fraction function $p_a(t)$. That is, the following Lemma holds.

*Lemma 1:* For any $a > 0$ and any real number $x_i, x_j$, the following inequalities hold:

$$\begin{aligned} p_a(x_i + x_j) &\le p_a(|x_i| + |x_j|) \\ &\le p_a(x_i) + p_a(x_j) \\ &\le 2p_a\left(\frac{|x_i| + |x_j|}{2}\right). \end{aligned} \tag{7}$$

In addition, because

$$p_a(t) \ge 0 \quad \text{and } p_a(t) = 0 \Leftrightarrow t = 0$$

the fraction function $p_a(t)$ acts almost like a norm. However, it lacks homogeneity, i.e., $p_a(ct) = |c| p_a(t)$, in general. In fact, the following Lemma is obvious.

*Lemma 2:* For the fraction function $p_a(t)$, we have

$$p_a(ct) < |c| p_a(t) \quad \text{if } |c| > 1$$

and

$$p_a(ct) \ge |c| p_a(t) \quad \text{if } |c| \le 1.$$

## III. Equivalence of the Minimization Problem $(FP_a)$ and $(P_0)$

We shall establish the equivalence (Theorems 1 and 2) of the minimization problem $(FP_a)$ and $(P_0)$ in this section.

*Lemma 3:* Let $x^*$ be the optimal solution to $(FP_a)$. Then, the columns in matrix $A$ corresponding to the support of vector $x^*$ are linearly independent, and hence, $\|x^*\|_0 = k \le m$.

*Proof:* See the Appendix.

We denote by $\Delta$ the set of solutions to $Ax = b$ with $\|x\|_0 \le m$, that is

$$\Delta = \{x | Ax = b \text{ and } \|x\|_0 \le m\}.$$

Obviously, the cardinality of $\Delta$ is finite.

*Theorem 1:* There exists some constant $\hat{a} > 0$, such that the optimal solution to the minimization problem $(FP_{\hat{a}})$ also solves the minimization problem $(P_0)$.

*Proof:* Let $\{a_i | i = 1, 2, \ldots, n, \ldots\}$ be the increasing infinite sequence with $\lim_{i\to\infty} a_i = \infty$ and $a_1 = 1$. For each $a_i$, by Lemma 3, the optimal solution $\hat{x}_{a_i}$ to $(FP_{a_i})$ belongs to the set $\Delta$. Since $\Delta$ is a finite set, one element, named $\hat{x}$, will repeatedly solve $(FP_{a_{i_k}})$ for some subsequence $\{a_{i_k} | k = 1, 2, \ldots\}$ of $\{a_i\}$. For any $a_{i_k} \geq a_{i_1}$, $x$ with $Ax = b$, we have

$$P_{a_{i_k}}(\hat{x}) = P_{a_{i_k}}(\hat{x}_{a_{i_k}}) = \min_{Ax=b} P_{a_{i_k}}(x) \leq \|x\|_0.$$

Letting $k \to \infty$ in the equality earlier, we have

$$\lim_{k\to\infty} P_{a_{i_k}}(\hat{x}) \leq \|x\|_0$$

that is, $\|\hat{x}\|_0 \leq \|x\|_0$. Hence, $\hat{x}$ is the optimal solution to $(P_0)$, which means that the optimal solution $\hat{x}$ to $(FP_{a_{i_k}})$ solves $(P_0)$. In general, we denote the smallest number $a_{i_1}$ of the infinite subsequence $\{a_{i_k} | k = 1, 2, \ldots\}$ by $\hat{a}$. The proof is completed.

It should be pointed out that the constant $\hat{a}$ in Theorem 1 may be very small. Before the following theorem is proved, a constant $r(A, b)$ needs to be defined. Let

$$r(A, b) = \min_{z\in\Delta, z_i\neq 0, 1\leq i\leq n} |z_i|. \quad (8)$$

Clearly, the defined constant $r(A, b)$ is finite and positive due to the finiteness of $\Delta$.

*Theorem 2:* There exists a constant $a^* > 0$, such that whenever $a > a^*$, every optimal solution to $(FP_a)$ also solves $(P_0)$, where $a^*$ depends on $A$ and $b$.

*Proof:* Let $x^*$ be the optimal solution to $(FP_a)$ and $x^0$ be the optimal solution to $(P_0)$. By Lemma 3, we know that $x^* \in \Delta$.

Therefore, we have

$$\begin{aligned}\min_{Ax=b} \|x\|_0 &= \|x^0\|_0 \\ &\geq \sum_{i\in\mathrm{supp}(x^0)} \frac{a|x_i^0|}{1 + a|x_i^0|} \\ &\geq \sum_{i\in\mathrm{supp}(x^*)} \frac{a|x_i^*|}{1 + a|x_i^*|} \\ &\geq \|x^*\|_0 \frac{ar(A, b)}{1 + ar(A, b)}\end{aligned}$$

which implies that

$$\|x^*\|_0 \leq \left(1 + \frac{1}{ar}\right)\|x^0\|_0 = \left(1 + \frac{1}{ar}\right)\min_{Ax=b} \|x\|_0.$$

Because $\|x^*\|_0$ is an integer number, it follows from the aforementioned inequality that $\|x^*\|_0 = \min_{Ax=b} \|x\|_0$ (that is, $x^*$ solves $(P_0)$) when:

$$\left(1 + \frac{1}{ar(A, b)}\right)\min_{Ax=b} \|x\|_0 < \min_{Ax=b} \|x\|_0 + 1. \quad (9)$$

Obviously, inequality (9) is true whenever

$$a > \frac{\min_{Ax=b} \|x\|_0}{r(A, b)}. \quad (10)$$

Therefore, with $a^*$ denoting the right-hand side of inequality (10), we conclude that when $a > a^*$, every solution $x^*$ to $(FP_a)$ also solves $(P_0)$. The proof is thus completed.

Although both Theorems 1 and 2 describe the equivalence of the minimization problem $(FP_a)$ and $(P_0)$, the constant $\hat{a}$ in Theorem 1 may be very small, and the constant $a^*$ in Theorem 2 is generally very large. Therefore, Theorem 2 is mainly used in theoretical analysis; however, when conducting the experiments, we often choose a smaller $a$ value by Theorem 1. In Section VI, we take $a = 2$, and the $FP$ algorithm performs well.

## IV. Equivalence of the Minimization Problem $(FP_a)$ and $(FP_a^\lambda)$

In this section, we first discuss the properties of the optimal solution to the regularization problem $(FP_a^\lambda)$, including the first-order and second-order optimality conditions and the lower and upper bounds of the absolute value for its nonzero entries. Second, based on these lemmas, we demonstrate the equivalence of the regularization problem $(FP_a^\lambda)$ and the constrained problem $(FP_a)$.

Before we embark on the discussion, we should mention that the results derived in this section are worst case ones, implying that the kind of guarantees we obtain are overly pessimistic, as they are supposed to hold for all signals and for all possible supports of a given cardinality.

*Lemma 4:* Suppose that $x^*$ is the optimal solution to $(FP_a^\lambda)$. Then, the following statements hold.

1) If $\lambda > \|b\|_2^2$, then

$$\|x^*\|_\infty \leq \frac{\|b\|_2^2}{a(\lambda - \|b\|_2^2)}.$$

2) Let $B$ be the submatrix of $A$ corresponding to the support of vector $x^*$. Thus, the columns of $B$ are linearly independent, and hence, $\|x^*\|_0 \leq m$.

3) Denote by $\bar{\lambda}$ the constant

$$\|b\|_2^2 + \frac{\|A^Tb\|_\infty + \sqrt{\|A^Tb\|_\infty + 2a\|b\|_2^2\|A^Tb\|_\infty}}{a}.$$

Then, for all $\lambda \geq \bar{\lambda}$, $x^* = 0$.

*Proof:* See the Appendix.

*Lemma 5 (First-Order Optimality Condition):* Let $x^*$ be the solution to $(FP_a^\lambda)$. Then, the following statements hold.

1) For any $h \in \mathbb{R}^n$ with $\mathrm{supp}(h) \subseteq \mathrm{supp}(x^*)$

$$2\langle b - Ax^*, Ah\rangle = \lambda \sum_{i\in\mathrm{supp}(x^*)} \frac{ah_i \mathrm{sgn}(x_i^*)}{(1 + a|x_i^*|)^2}. \quad (11)$$

2) For any $h \in \mathbb{R}^n$ with $\mathrm{supp}(h) \subseteq \mathrm{Csupp}(x^*)$

$$|\langle b - Ax^*, Ah\rangle| \leq \sqrt{\|Ah\|_2^2\lambda\|h\|_0} \quad (12)$$

where $\mathrm{Csupp}(x^*)$ is the complementary of $\mathrm{supp}(x^*)$.

*Proof:* See the Appendix.

Choosing $h$ as the $i$th base vector $e_i$ for each $i = 1, 2, \ldots, n$ in (11) and (12), respectively, we can derive the following corollary.

*Corollary 1:* Suppose that $x^*$ is the solution to $(FP_a^\lambda)$. Then, for $i \in \text{supp}(x^*)$

$$2(A^T(b - Ax^*))_i = \lambda \frac{a\text{sgn}(x_i^*)}{(1 + a|x_i^*|)^2} \tag{13}$$

and for $i \in \text{Csupp}(x^*)$

$$|(A^T(b - Ax^*))_i| \le \sqrt{\lambda \|a_i\|_2^2}.$$

Furthermore, letting $\lambda > \|b\|_2^2$ and replacing $|x_i^*|$ with $\|x^*\|_\infty$ in (13), we have

$$\|A^T(b - Ax^*)\|_2^2 \ge \frac{a^2(\lambda - \|b\|_2^2)^4}{4\lambda^2}. \tag{14}$$

Following the above-adopted analysis, we can further establish the following optimality condition.

*Lemma 6 (Second-Order Optimality Condition):* Every solution $x^*$ to $(FP_a^\lambda)$ satisfies the following conditions.

1) For all $h \in \mathbb{R}^n$ with $\text{supp}(h) \subseteq \text{supp}(x^*)$

$$\|Ah\|_2^2 \ge \lambda \sum_{i \in \text{supp}(x^*)} \frac{2a^2 h_i^2}{(1 + a|x_i^*|)^3}. \tag{15}$$

2) Moreover, it holds for all $i \in \text{supp}(x^*)$ that

$$|x_i^*| \ge \frac{\sqrt{\lambda}}{\|a_i\|_2} - \frac{1}{a}. \tag{16}$$

*Proof:* See the Appendix.

In the following, we discuss the equivalence of the regularization problem $(FP_a^\lambda)$ and the constrained problem $(FP_a)$. We denote by $\sigma_{\min}$ the minimal one of all the smallest singular values of $A_s$, where $A_s$ is an arbitrary submatrix of $A$ with full column rank. That is

$$\sigma_{\min} = \min\{\sigma_s | \sigma_s \text{ is the smallest singular value of } A_s, \text{ where } A_s \text{ is an arbitrary submatrix of } A \text{ with full column rank}\}. \tag{17}$$

Clearly, $\sigma_{\min} > 0$.

*Theorem 3:* If there exists a constant $\lambda \in (\|b\|_2^2, \bar{\lambda})$, such that

$$\frac{4m\|A\|_2^4}{\lambda a^2} \left( \frac{\lambda}{\lambda - \|b\|_2^2} \right)^4 < \sigma_{\min} \tag{18}$$

then the optimal solution to $(FP_a^\lambda)$ solves $(FP_a)$, where $\bar{\lambda}$ is defined in Lemma 4.

*Proof:* See the Appendix.

Moreover, if the constant $a$ in Theorem 3 satisfies $a > a^*$ ($a^*$ is the one in Theorem 2), then we have the following corollary by Theorems 2 and 3.

*Corollary 2:* If the constant $a$ in $(FP_a^\lambda)$ satisfies $a > a^*$ and there exists a constant $\lambda \in (\|b\|_2^2, \bar{\lambda})$, such that (18) holds, then the optimal solution to $(FP_a^\lambda)$ also solves $(P_0)$, where $\bar{\lambda}$ is defined in Lemma 4.

Theorem 3 and Corollary 2 show that it is possible to obtain the exact solution to $(P_0)$ by solving the problem $(FP_a^\lambda)$. In Section V, we will discuss the algorithms to solve $(FP_a^\lambda)$.

## V. Thresholding Algorithms for the Regularization Problem $(FP_a^\lambda)$

In this section, we mainly derive the closed-form representation of the optimal solution to the regularization problem $(FP_a^\lambda)$ (see Theorem 4), which underlies the algorithm to be proposed. Before giving the closed-form representation of the optimal solution, we need to prove Lemmas 8–10.

Let us define three parameters $t_1^*, t_2^*$, and $t_3^*$ for our following derivation, where:

$$t_1^* = \frac{\sqrt[3]{\frac{27}{8}\lambda a^2} - 1}{a} \quad t_2^* = \frac{\lambda}{2}a \quad \text{and } t_3^* = \sqrt{\lambda} - \frac{1}{2a}.$$

Obviously, we have the following lemmas.

*Lemma 7:* For any positive parameters $\lambda$ and $a$, $t_1^* \le t_3^* \le t_2^*$ holds. Furthermore, they are equal to $1/2a$ when $\lambda = 1/a^2$.

*Lemma 8:* For any given $t$, the two polynomials of $x$ defined in the following satisfy the following conditions.

1) If $t > t_1^*$, then the polynomial

$$2x(ax + 1)^2 - 2t(ax + 1)^2 + \lambda a = 0 \tag{19}$$

has three different real roots, and the largest root $x_0$ is obtained by $x_0 = g_\lambda(t)$, where

$$g_\lambda(t) = \text{sgn}(t) \left( \frac{\frac{1+at}{3}\left(1 + 2\cos\left(\frac{\phi(t)}{3} - \frac{\pi}{3}\right)\right) - 1}{a} \right)$$

$$\phi(t) = \arccos\left( \frac{27\lambda a^2}{4(1 + a|t|)^3} - 1 \right).$$

Clearly, $|g_\lambda(t)| \le |t|$.

2) If $t < -t_1^*$, then

$$2x(1 - ax)^2 - 2t(1 - ax)^2 - \lambda a = 0 \tag{20}$$

has three different real roots, and the smallest root $x_0$ is obtained by $x_0 = g_\lambda(t)$.

*Proof:* See the Appendix.

We define a function of $y$ as

$$f_\lambda(y) = (y - x)^2 + \lambda p_a(|y|).$$

*Lemma 9:* The optimal solution to $\min_{y \in \mathbb{R}} f_\lambda(y)$ is the threshold function defined as

$$y^* = \begin{cases} 0, & |x| \le t \\ g_\lambda(x), & |x| > t \end{cases} \tag{21}$$

where $g_\lambda(x)$ is the one in Lemma 8, and the parameter $t$ satisfies

$$t = \begin{cases} t_2^*, & \text{if } \lambda \le \frac{1}{a^2} \\ t_3^*, & \text{if } \lambda > \frac{1}{a^2}. \end{cases}$$

*Proof:* See the Appendix.

For any $\lambda,\ \mu \in (0, +\infty)$ and $z \in \mathbb{R}^n$, let

$$C_\lambda(x) = \|Ax - b\|_2^2 + \lambda P_a(x) \tag{22}$$

$$C_\mu(x, z) = \mu\big(C_\lambda(x) - \|Ax - Az\|_2^2\big) + \|x - z\|_2^2 \tag{23}$$

and

$$B_\mu(x) = x + \mu A^T(b - Ax). \tag{24}$$

*Lemma 10:* For any fixed parameter $\mu$, $a$, $\lambda$, and $z$, if $x^s = (x_1^s, x_2^s, \dots, x_n^s)^T$ is a local minimizer to $C_\mu(x, z)$, then

$$x_i^s = 0 \Leftrightarrow |(B_\mu(z))_i| \le t^*$$

and

$$x_i^s = g_{\lambda\mu}((B_\mu(z))_i) \Leftrightarrow |(B_\mu(z))_i| > t^*$$

where

$$t^* = \begin{cases} \dfrac{\lambda\mu a}{2}, & \text{if } \lambda \le \dfrac{1}{\mu a^2} \\ \sqrt{\lambda\mu} - \dfrac{1}{2a}, & \text{if } \lambda > \dfrac{1}{\mu a^2} \end{cases}$$

and $g_{\lambda\mu}(\cdot)$ is obtained by replacing $\lambda$ with $\lambda\mu$ in $g_\lambda(\cdot)$.

*Proof:* We first notice that $C_\mu(x, z)$ can be rewritten as

$$\begin{aligned} C_\mu(x, z) &= \|x - ((I - \mu A^T A)z + \mu A^T b)\|_2^2 \\ &\quad + \lambda\mu P_a(x) + \mu\|b\|_2^2 + \|z\|_2^2 - \mu\|Az\|_2^2 \\ &\quad - \|(I - \mu A^T A)z + \mu A^T b)\|_2^2 \\ &= \sum_{i=1}^n (x_i - (B_\mu(z))_i)^2 + \lambda\mu \sum_{i=1}^n p_a(|x_i|) \\ &\quad + \mu\|b\|_2^2 + \|z\|_2^2 - \mu\|Az\|_2^2 \\ &\quad - \|(I - \mu A^T A)z + \mu A^T b)\|_2^2 \end{aligned}$$

which implies that minimizing $C_\mu(x, z)$ for any fixed $\mu$, $\lambda$, and $z$ is equivalent to

$$\min_{x \in \mathbb{R}^n} \sum_{i=1}^n (x_i - (B_\mu(z))_i)^2 + \lambda\mu \sum_{i=1}^n p_a(|x_i|).$$

Therefore, $x^s = (x_1^s, x_2^s, \dots x_n^s)^T$ is a local minimizer of $C_\mu(x, z)$ if and only if, for any $i, x_i^s$ solves the problem

$$\min_{x \in \mathbb{R}^n} (x_i - (B_\mu(z))_i)^2 + \lambda\mu p_a(|x_i|).$$

Therefore, the proof is completed by Lemma 9.

Now, we show that the optimal solution to the problem $(FP_a^\lambda)$ can be expressed as a thresholding operation.

*Theorem 4:* If $x^* = (x_1^*, x_2^*, \dots x_n^*)^T$ is an optimal solution to $(FP_a^\lambda)$, $a$ and $\lambda$ are positive, and parameter $\mu$ satisfies $0 < \mu < \|A\|_2^{-2}$, then the optimal solution $x^*$ is

$$x_i^* = \begin{cases} g_{\lambda\mu}((B_\mu(x^*))_i), & |(B_\mu(x^*))_i| > t^* \\ 0, & \text{otherwise} \end{cases} \tag{25}$$

where parameter $t^*$ is defined in Lemma 10.

*Proof:* The condition $0 < \mu < \|A\|_2^{-2}$ implies that

$$\begin{aligned} C_\mu(x, x^*) &= \mu\big(\|b - Ax\|_2^2 + \lambda P_a(x)\big) \\ &\quad + \big(-\mu\|Ax - Ax^*\|_2^2 + \|x - x^*\|_2^2\big) \\ &\ge \mu\big(\|b - Ax\|_2^2 + \lambda P_a(x)\big) \\ &\ge C_\mu(x^*, x^*) \end{aligned}$$

for any $x \in \mathbb{R}^n$, which shows that $x^*$ is a local minimizer of $C_\mu(x, x^*)$ as long as $x^*$ is a solution to $(FP_a^\lambda)$. Following directly from Lemmas 9 and 10, we finish the proof.

In the following, we present an iterative thresholding algorithm for performing the regularization problem $(FP_a^\lambda)$ based on the previous theoretical analysis.

With the thresholding representation (25), a thresholding algorithm for the regularization problem $(FP_a^\lambda)$ can be naturally defined as

$$x_i^{n+1} = \begin{cases} g_{\lambda\mu}((B_\mu(x^n))_i), & |(B_\mu(x^n))_i| > t^* \\ 0, & \text{otherwise} \end{cases} \tag{26}$$

where $g_{\lambda\mu}(\cdot)$ is the thresholding operator defined in Lemma 10. We call this method the iterative $FP$ thresholding algorithm, or briefly, the $FP$ algorithm.

It is known that the quality of the solutions of a regularization problem depends seriously on the setting of the regularization parameter $\lambda$. However, the selection of proper regularization parameters is a very hard problem. In most general cases, a "trial and error" method, for example, the cross-validation method, is still an accepted, or even unique, choice. Nevertheless, when some prior information is known for a problem, it is realistic to set the regularization parameter more reasonably and intelligently.

To make it clear, let us suppose that the solutions to the regularization problem $FP_a^\lambda$ are of $k$-sparsity. Thus, we are required to solve the regularization problem $FP_a^\lambda$ restricted to the subregion $\Gamma_k = \{x = (x_1, x_2, \dots, x_n) | \text{supp}(x) = k\}$ of $\mathbb{R}^n$. Assume that $x^*$ is the solution to the regularization problem $FP_a^\lambda$, and without loss of generality, $|B_\mu(x^*)|_1 \ge |B_\mu(x^*)|_2 \ge \cdots \ge |B_\mu(x^*)|_n$. Then, by Theorem 4, the following inequalities hold:

$$|B_\mu(x^*)|_i > t^* \Leftrightarrow i \in 1, 2, \dots, k$$
$$|B_\mu(x^*)|_j \le t^* \Leftrightarrow j \in k+1, k+2, \dots, n$$

where $t^*$ is our threshold value that was defined earlier. According to $t_3^* \le t_2^*$, we have

$$\begin{cases} |B_\mu(x^*)|_k > t^* \ge \sqrt{\lambda\mu} - \dfrac{1}{2a} \\ |B_\mu(x^*)|_{k+1} \le t^* \le \dfrac{\lambda\mu a}{2} \end{cases} \tag{27}$$

which implies

$$\frac{2|B_\mu(x^*)|_{k+1}}{a\mu} \le \lambda < \frac{(2a|B_\mu(x^*)|_k + 1)^2}{4a^2\mu}. \tag{28}$$

For convenience, we denote by $\lambda_1$ and $\lambda_2$ the left- and right-hand side of the above-mentioned inequality, respectively.

The aforementioned estimate helps to set an optimal regularization parameter. A choice of $\lambda$ is

$$\lambda = \begin{cases} \lambda_1, & \lambda_1 \le \dfrac{1}{a^2\mu} \\ (1 - \epsilon)\lambda_2, & \lambda_1 > \dfrac{1}{a^2\mu} \end{cases}$$

where $\epsilon$ is a small positive number, such as 0.1, 0.01, or 0.001.

**Algorithm 1** Iterative FP Thresholding Algorithm—Scheme 1

Initialization: Choose $x^0$, $\mu = \mu_0 = \frac{1-\varepsilon}{\|A\|^2}$ and $a$;
**while** not converged **do**
  $z^n := B_\mu(x^n) = x^n + \mu A^T(b - Ax^n)$;
  $\lambda = \lambda_0$;
  if $\lambda \le \frac{1}{a^2\mu}$, then
    $t^* = \frac{\lambda\mu a}{2}$
    for $i = 1 : length(x)$
      $1.|z_i^n| > t^*, then \quad x_i^{n+1} = g_{\lambda\mu}(z_i^n)$
      $2.|z_i^n| \le t^*, then \quad x_i^{n+1} = 0$
  else
    $t^* = \sqrt{\lambda\mu} - \frac{1}{2a}$
    for $i = 1 : length(x)$
      $1.|z_i^n| > t^*, then \quad x_i^{n+1} = g_{\lambda\mu}(z_i^n)$
      $2.|z_i^n| \le t^*, then \quad x_i^{n+1} = 0$
  end
  $n \to n + 1$
**end while**
**return** $x^{n+1}$

In practice, we approximate $x^*$ by $x^n$ in (28). Suppose we take

$$\lambda^* = \begin{cases} \lambda_1 = \dfrac{2|B_\mu(x^n)|_{k+1}}{a\mu}, & \lambda_1 \le \dfrac{1}{a^2\mu} \\ (1-\epsilon)\lambda_2 = (1-\epsilon)\dfrac{(2a|B_\mu(x^n)|_k + 1)^2}{4a^2\mu}, & \lambda_1 > \dfrac{1}{a^2\mu} \end{cases} \tag{29}$$

in applications. When doing so, an iteration algorithm will be adaptive and free from the choice of the regularization parameter. Note that (29) is valid for any $\mu$ satisfying $0 < \mu < \|A\|_2^{-2}$. In general, we can take $\mu = \mu_0 = 1 - \varepsilon/\|A\|_2^2$ with any small $\varepsilon \in (0, 1)$ below.

Incorporated with different parameter-setting strategies, (26) defines different implementation schemes of the $FP$ algorithm. For example, we can have the following.

*Scheme 1:* $\mu = \mu_0$; $\lambda_n = \lambda_0 \in (\|b\|_2^2, \bar{\lambda})$ and $a = a_0$.

*Scheme 2:* $\mu = \mu_0$; $\lambda_n = \lambda^*$ defined in (29) and $a = a_0$.

Importantly, it should be noted that the threshold value $t^* = \lambda\mu a/2$ when the parameter $\lambda_n = \lambda_1$, and the threshold value $t^* = (\lambda\mu)^{1/2} - 1/2a$ when the parameter $\lambda_n = (1-\epsilon)\lambda_2$ in Scheme 2. Our analysis leads to the algorithm in Algorithms 1 and 2.

In Algorithm 2, $\lambda$ is tuned with Scheme 2; the others are the same as Algorithm 1.

At the end of the section, we discuss the convergence of the $FP$ algorithm to a stationary point of the iteration (26) under some certain conditions.

*Theorem 5:* Let $\{x^k\}$ be the sequence generated by the $FP$ algorithm with $0 < \mu < \|A\|_2^{-2}$. Then, the following statements hold.

1) The sequence $C_\lambda(x^k) = \|Ax^k - b\|_2^2 + \lambda P_a(x^k)$ is decreasing.
2) $\{x^k\}$ is asymptotically regular, i.e., $\lim_{k\to\infty} \|x^{k+1} - x^k\|_2 = 0$.
3) $\{x^k\}$ converges to a stationary point of the iteration (26).

*Proof:* The proof is similar to [25, Proof of Lemma 2], so it is omitted here.

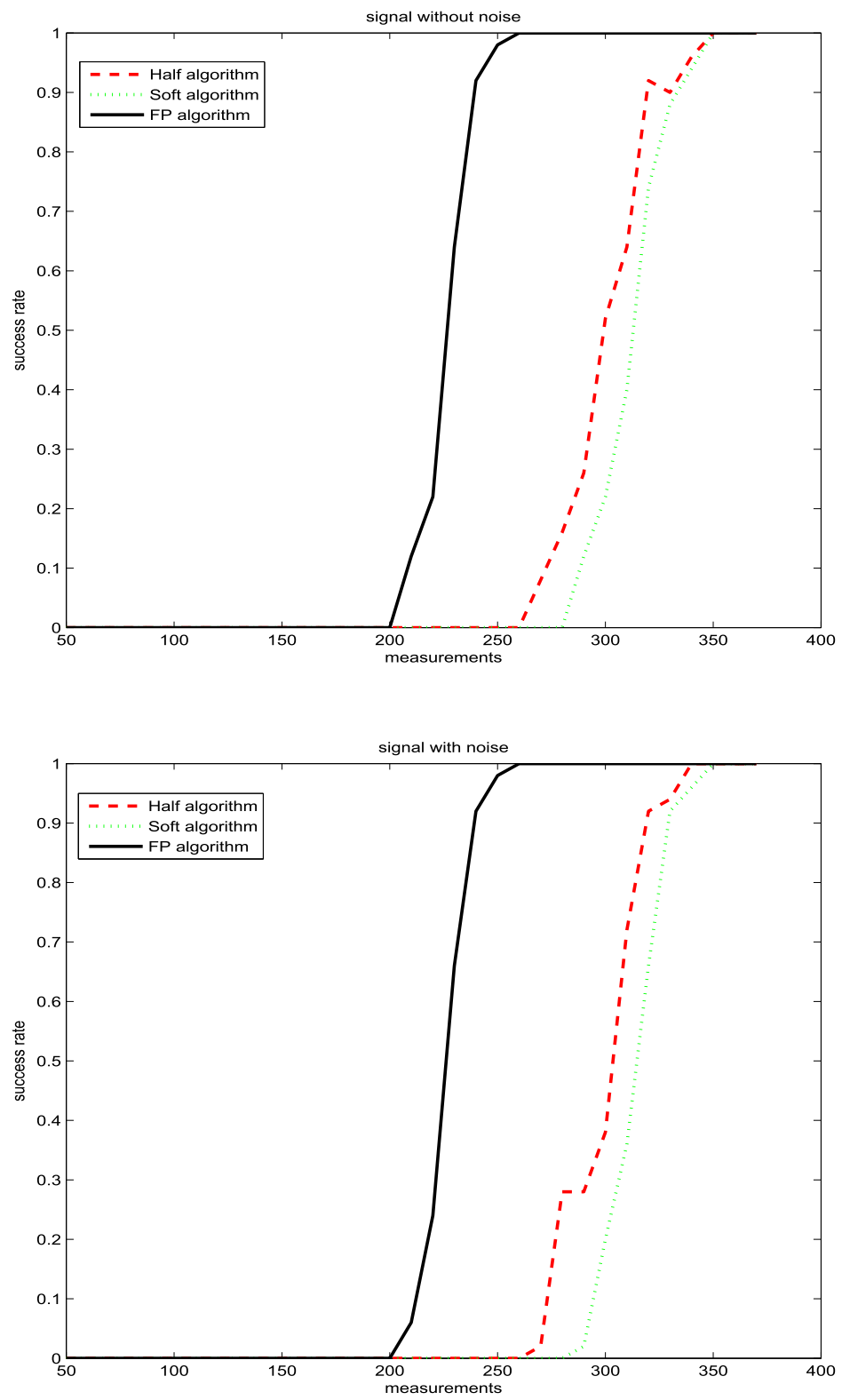

Fig. 1. How few measurements (samples) by three algorithms are required to exactly recover a given signal $x_0$.

## VI. Experimental Results

In this section, we carry out a series of simulations to demonstrate the performance of the $FP$ algorithm. All the simulations here are conducted by applying our algorithm (see Scheme 2) to a typical compressed sensing problem, i.e., signal recovery. In the experiments, the soft algorithm, the half algorithm, and the $FP$ algorithm are simulated from three aspects. For each experiment, we repeatedly perform 100 tests, present the average results, and take $a = 2$.

The simulations are all conducted on a personal computer (3.60 GHz, 4-GB RAM) with MATLAB 8.0 programming platform (R2012b).

The first simulation tests how few measurements (samples) are required to exactly recover a given signal $x_0$ using the three algorithms. It is obvious that the fewer the measurements that are used by an algorithm, the better it is. Consider a real-valued $N$-length ($N = 512$) signal $x_0$ without noise, which is randomly generated under a Gaussian distribution of zero mean and unit variance, $N(0, 1)$, and its sparsity is fixed at $k = 100$. The simulations then aim to recover $x \in \mathbb{R}^{512}$ through $M$ measurements determined by the measurement matrix $A_{M\times512}$, where $A_{M\times512}$ is a random matrix with entries independently drawn at random from a Gaussian distribution of zero mean and unit variance, $N(0, 1)$, and $M$ ranges from 50 to 370. The three algorithms are applied with a variable

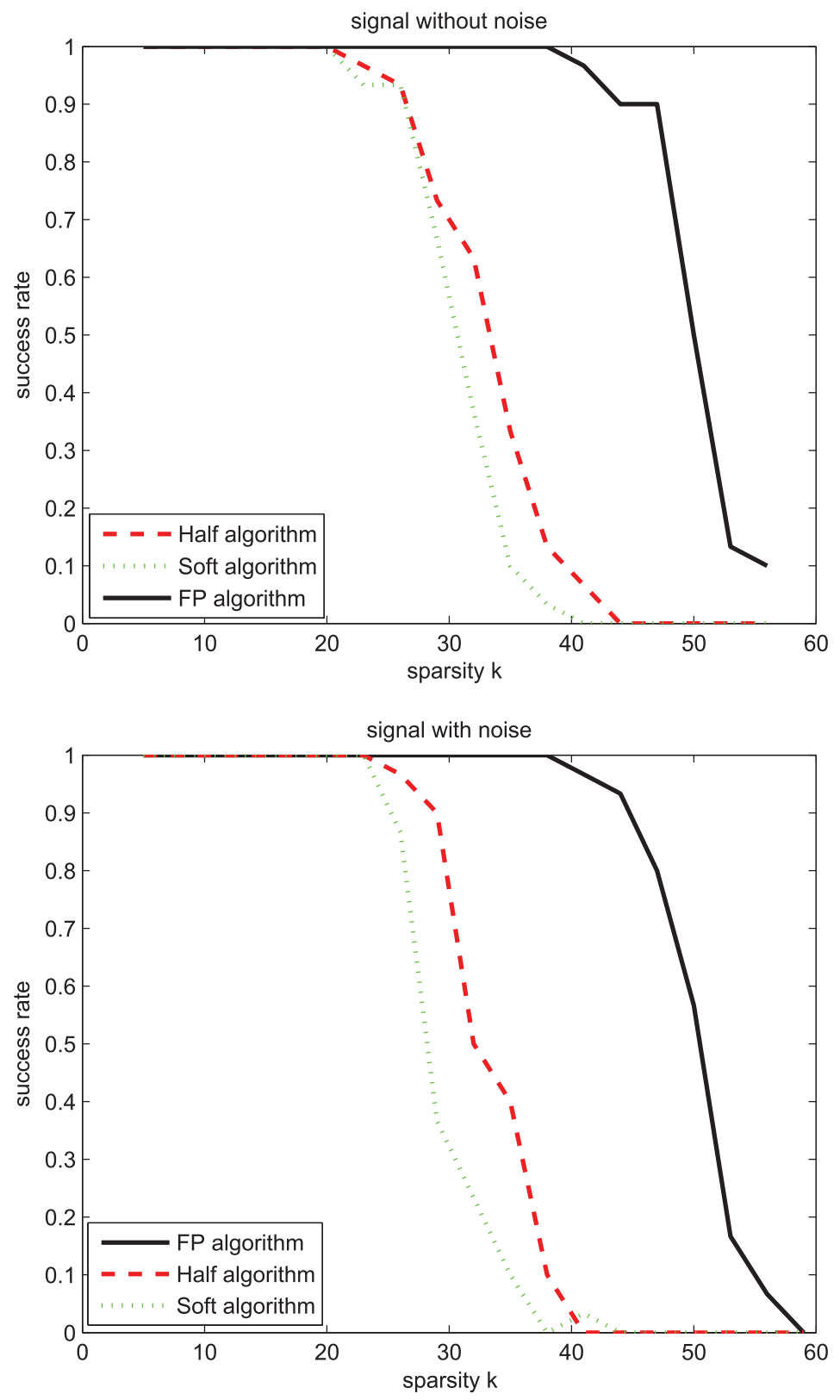

Fig. 2. Success rate of the three algorithms in the recovery of a signal $x_0$ with a different cardinality for a given measurement matrix $A$.

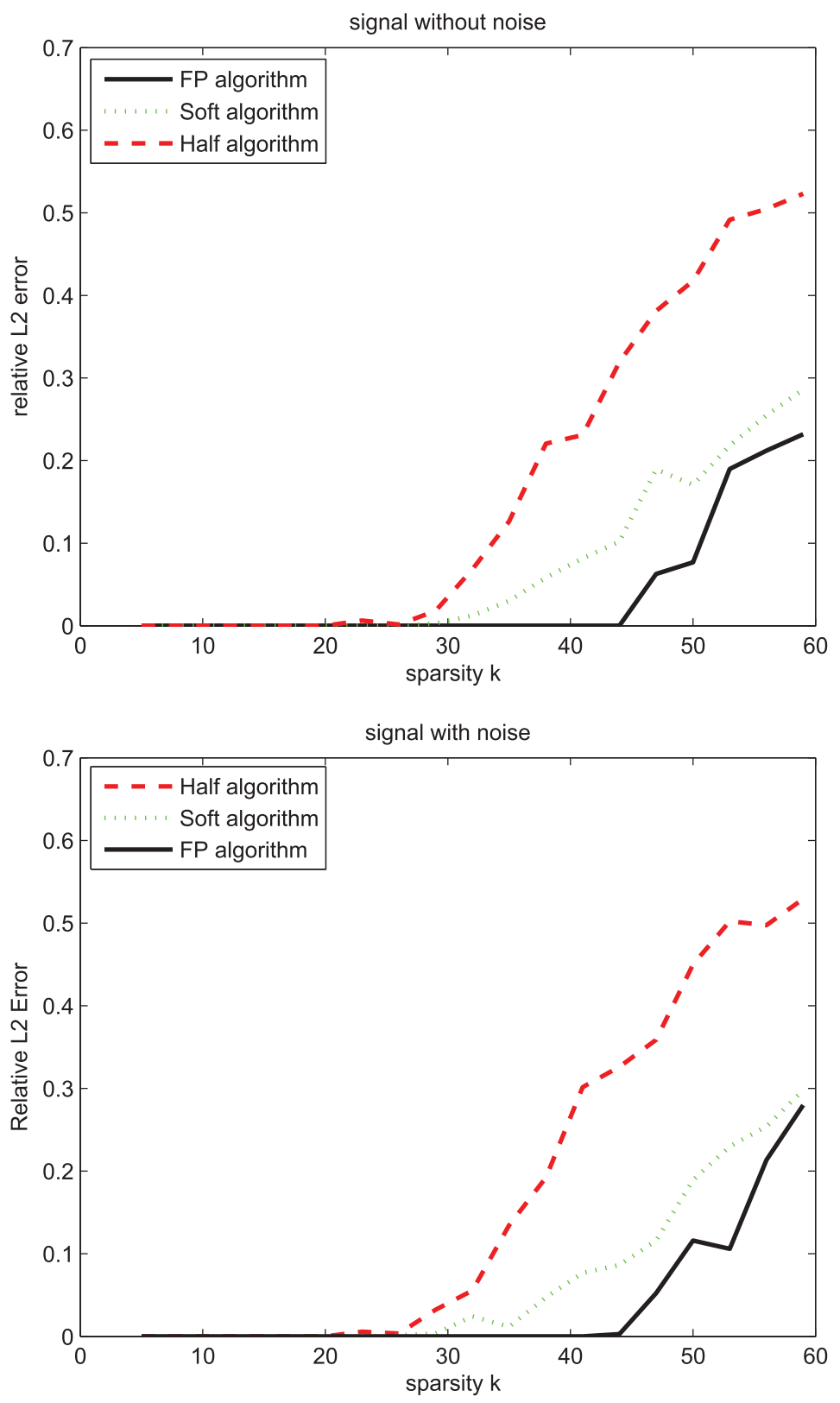

Fig. 3. Relative $\ell_2$-error between the solution $\hat{x}$ and the given signal $x_0$.

number $M$ of measurements. The results of the simulations are shown in Fig. 1.

Turning to the noisy case, we use the same signal $x_0$ but with noise, for example, with the white noise $\varepsilon \in N(0, \sigma^2)$ ($\sigma = 0.1$). Such a noise signal is designed to simulate a real measurement in which noise is inevitably involved. Our simulations aim to assess the capability of all three algorithms in recovering the signal from a noisy circumstance and with fewer samplings. The results of the simulations are shown in Fig. 1.

From Fig. 1, we can see that the three algorithms can accurately recover the signal $x_0$ when $M \geq 350$, and when the measurements are deduced to 260, there is no other algorithm except the $FP$ algorithm that can accurately recover the signal $x_0$. The simulation results show that the $FP$ algorithm requires the least number of samplings among the three algorithms. The graph presented in Fig. 1 shows that the $FP$ algorithm in recovering the signal from a noisy circumstance also requires the least number of samplings among the three algorithms.

The second simulation tests the success rate of the three algorithms in the recovery of a signal with a different cardinality for a given measurement matrix $A$. Consider a random matrix $A$ of size $128 \times 512$, with entries independently drawn by random from a Gaussian distribution of zero mean and unit variance, $N(0, 1)$. By randomly generating such sufficiently sparse vectors $x_0$ [choosing the nonzero locations uniformly over the support randomly, and their values from $N(0, 1)$], we generate vector $b$. This way, we know the sparsest solution to $Ax_0 = b$, and we are able to compare it to the algorithmic results. The success is measured by computing $(\|\hat{x} - x_0\|_2^2 / \|x_0\|_2^2)$ and checking that is below a negligible value (in our experiments, it is set to $1e^{-5}$), to indicate a perfect recovery of the original sparse vector $x_0$.

Turning to the noisy case, we use the same matrix $A$ and generate a random vector $x_0$ with a pre-specified cardinality of nonzeros. We compute $b = Ax_0 + \varepsilon$, where $\varepsilon \in N(0, \sigma^2)$ ($\sigma = 0.1$). Thus, the original vector $x_0$ is a feasible solution and close to the optimal solution. Due to the presence of noise, it becomes harder to accurately recover the original signal $x_0$. Therefore, we tune down the requirement for a success to the relative error $(\|x^* - x_0\|_2^2 / \|x_0\|_2^2 \leq 10^{-5})$.

The graphs presented in Fig. 2 show the success rate of the soft algorithm, the half algorithm, and the $FP$ algorithm in recovering the true (sparsest) solution. From Fig. 2, we can see that the $FP$ algorithm can exactly recover the ideal signal until $k$ is approximately 39, and the soft algorithm and half algorithm's counterpart is approximately slightly higher than 21. The results in the noisy state are consistent with the above-said one. As we can see, the $FP$ algorithm again has the best performance, with the half algorithm as the second best.

Next, we consider the relative $\ell_2$-error between the solution $\hat{x}$ and the given signal $x_0$. The $\ell_2$-error is computed as the ratio $(\|\hat{x} - x_0\|_2^2 / \|x_0\|_2^2)$, indicating the $\ell_2$-proximity between the two solutions, and we measured this distance as relative to the energy in the true solution. The results of the simulations are shown in Fig. 3.

From Fig. 3, we can see that the $FP$ algorithm always has the smallest relative $\ell_2$-error value.

## VII. Conclusion

As is well known, $(P_0)$ is combinatorial and NP-hard in general. Therefore, it is important to choose suitable substitution models for $\ell_0$ minimization. In this paper, we take the fraction function as the substitution for the $\ell_0$-norm and study the fraction function minimization in terms of theory, algorithms, and computation. In the beginning, we discuss the equivalence of $\ell_0$ minimization and fraction function minimization. In particular, we consider the regularization model $(FP_a^\lambda)$ and prove that under certain conditions, the optimal solution to $(FP_a^\lambda)$ also solves $(P_0)$.

The above-mentioned conclusions demonstrate that we can obtain the exact solution to $(P_0)$ by solving the regularization model $(FP_a^\lambda)$. Hence, it is necessary to study the algorithm for solving the regularization problem $(FP_a^\lambda)$. We derive the closed-form representation of the optimal solution to the regularization problem $(FP_a^\lambda)$. Based on this representation, we calculate its analytic expression and propose an iterative $FP$ thresholding algorithm to solve the regularization problem $(FP_a^\lambda)$.

We also provide a series of experiments to assess the performance of the $FP$ algorithm, and the experimental results show that the $FP$ algorithm performs well, although the solution to $(FP_a^\lambda)$ by the $FP$ algorithm is a stationary point and may not be the minimal solution in theory because the fraction function is non-convex. Therefore, in the future, we will investigate the sufficient conditions, such that the solution to $(FP_a^\lambda)$ by the $FP$ algorithm is the minimal one.

## Appendix

### Proofs of Theorems and Lemmas

*Proof of Lemma 3*: Let us assume that the $k$-sparsity vector $x^*$ is the optimal solution to $(FP_a)$ and the $k$ columns combined linearly by $x^*$ are linearly dependent. Without loss of generality, we assume that $x^* = (x_1^*, x_2^*, \ldots, x_k^*, 0, 0, \ldots, 0)$. Then, there exists a non-trivial vector $h^*$ that combines these columns to zero (i.e., the support of $h^*$ is contained within the support of $x^*$), $Ah^* = 0$. It is obvious that $A(x^*+h^*) = b$ and $A(x^* - h^*) = b$. Without loss of generality, we assume that $\max_{1\le j\le k} h_j^* \le \min_{1\le j\le k} x_j^*$. Hence, for every $j$, $x_j^* + h_j^*$, $x_j^* - h_j^*$, and $x_j^*$ have the same sign. Because the function $f(t) = (at/1 + at)(t > 0)$ is strictly concave, we have, for every $j$

$$\frac{a|x_j^* + h_j^*|}{1 + a|x_j^* + h_j^*|} + \frac{a|x_j^* - h_j^*|}{1 + a|x_j^* - h_j^*|} < 2\frac{a|x_j^*|}{1 + a|x_j^*|}.$$

Furthermore

$$\sum_{j=1}^k \frac{a|x_j^* + h_j^*|}{1 + a|x_j^* + h_j^*|} + \sum_{j=1}^k \frac{a|x_j^* - h_j^*|}{1 + a|x_j^* - h_j^*|} < 2\sum_{j=1}^k \frac{a|x_j^*|}{1 + a|x_j^*|}$$

which implies that

$$\sum_{j=1}^k \frac{a|x_j^* + h_j^*|}{1 + a|x_j^* + h_j^*|} < \sum_{j=1}^k \frac{a|x_j^*|}{1 + a|x_j^*|}$$

or

$$\sum_{j=1}^k \frac{a|x_j^* - h_j^*|}{1 + a|x_j^* - h_j^*|} < \sum_{j=1}^k \frac{a|x_j^*|}{1 + a|x_j^*|}.$$

That is

$$P_a(x^* + h^*) < P_a(x^*) \text{ or } P_a(x^* - h^*) < P_a(x^*).$$

This is a contraction.

*Proof of Lemma 4*:

1) Let $x^*$ be the optimal solution to $(FP_a^\lambda)$. Then, we have

$$f(x^*) = \|Ax^* - b\|_2^2 + \lambda P_a(x^*) \le f(0) = \|b\|_2^2.$$

Hence, $\lambda P_a(x^*) \le \|b\|_2^2$, which implies that

$$\frac{a\|x^*\|_\infty}{1 + a\|x^*\|_\infty} \le \frac{\|b\|_2^2}{\lambda}.$$

If $\lambda > \|b\|_2^2$, then

$$\|x^*\|_\infty \le \frac{\|b\|_2^2}{a(\lambda - \|b\|_2^2)}.$$

2) Let $B$ be the submatrix of $A$ corresponding to the support of vector $x^*$. By inequality (15) in Lemma 6, for any $y \ne 0$

$$\|By\|_2^2 \ge \lambda \sum_{i=1}^k \frac{2a^2 y_i^2}{(1 + a|x_i^*|)^3} > 0$$

which implies that the matrix $B^T B$ is positive definite. Thus, the columns of $B$ are linearly independent, and hence, $\|x^*\|_0 \le m$.

3) Suppose that $x^* \ne 0$ and $\|x^*\|_0 = k$. Without loss of generality, we assume

$$x^* = (x_1^*, \ldots, x_k^*, 0, \ldots, 0)^T.$$

Let $z^* = (x_1^*, \ldots, x_k^*)^T$ and $B \in \mathbb{R}^{m\times k}$ be the submatrix of $A$, whose columns are the first $k$ columns of $A$.

We define a function $g : \mathbb{R}^k \to \mathbb{R}$ by

$$g(z) = \|Bz - b\|_2^2 + \lambda P_a(z).$$

Then

$$\begin{aligned} f(x^*) &= \|Ax^* - b\|_2^2 + \lambda P_a(x^*) \\ &= \|Bz^* - b\|_2^2 + \lambda P_a(z^*) = g(z^*). \end{aligned}$$

Since $|z_i^*| > 0, i = 1, 2, \ldots, k,$ $g$ is continuously differentiable at $z^*$, and moreover, in a neighborhood of $x^*$

$$\begin{aligned} g(z^*) &= f(x^*) \le \min\{f(x)|x_i = 0, i = k+1, \ldots, n\} \\ &= \min\{g(z)|z \in \mathbb{R}^k\} \end{aligned}$$

which implies that $z^*$ is a local minimizer of the function $g$. Hence, the first-order necessary condition for

$$\min_{z\in\mathbb{R}^k} g(z)$$

at $z^*$ gives

$$2B^T(Bz^* - b) + \text{diag}(\text{sgn}(z))\frac{\lambda a}{(1 + a|z|)^2} = 0$$

where sgn(·) is the sign function. Multiplying by $z^{*T}$ both the sides of the above-mentioned equality yields

$$2z^{*T}B^TBz^* - 2z^{*T}B^Tb + z^{*T}\text{diag}(\text{sgn}(z))\frac{\lambda a}{(1+a|z|)^2} = 0.$$

Because the columns of $B$ are linearly independent, $B^TB$ is positive definite, and hence

$$-2z^{*T}B^Tb + z^{*T}\text{diag}(\text{sgn}(z^*))\frac{\lambda a}{(1+a|z^*|)^2} < 0.$$

Equivalently

$$\sum_{i=1}^{k}\left(\frac{\lambda a|z_i^*|}{\left(1+a|z_i^*|\right)^2} - 2(B^Tb)_iz_i^*\right) < 0. \tag{30}$$

Since

$$\lambda \geq \|b\|_2^2 + \frac{\|A^Tb\|_\infty + \sqrt{\|A^Tb\|_\infty + 2a\|b\|_2^2\|A^Tb\|_\infty}}{a}$$

we obtain

$$a\lambda^2 - 2\left(a\|b\|_2^2 + \|A^Tb\|_\infty\right)\lambda + a\|b\|_2^4 \geq 0 \tag{31}$$

which implies that

$$\frac{a\left(\lambda - \|b\|_2^2\right)^2}{\lambda} \geq 2\|A^Tb\|_\infty. \tag{32}$$

Together with

$$\frac{\lambda a}{\left(1+a|z_i^*|\right)^2} \geq \frac{a\left(\lambda - \|b\|_2^2\right)^2}{\lambda} \text{ and } 2|(B^Tb)_i| \leq 2\|A^Tb\|_\infty$$

we obtain that

$$\frac{\lambda a}{\left(1+a|z_i^*|\right)^2} - 2|(B^Tb)_i| \geq 0. \tag{33}$$

Hence, for any $i \in \text{supp}(z^*)$

$$\frac{\lambda a|z_i^*|}{\left(1+a|z_i^*|\right)^2} - 2(B^Tb)_iz_i^* \geq 0$$

which is a contraction with (30), as claimed.

*Proof of Lemma 5*: Let $x^*$ be any solution to $(FP_a^\lambda)$. Then, for all $t \in \mathbb{R}$ and $h \in \mathbb{R}^n$, the following inequality holds:

$$\|Ax^* - b\|_2^2 + \lambda P_a(x^*) \leq \|A(x^* + th) - b\|_2^2 + \lambda P_a(x^* + th).$$

Equivalently

$$t^2\|Ah\|_2^2 + 2t\langle Ax^* - b, Ah\rangle + \lambda(P_a(x^* + th) - P_a(x^*)) \geq 0. \tag{34}$$

1) If $\text{supp}(h) \subseteq \text{supp}(x^*)$, then for all $t \in \mathbb{R}$

$$P_a(x^* + th) - P_a(x^*) = \sum_{i\in\text{supp}(x^*)}\left(\frac{a|x_i^* + th_i|}{1+a|x_i^* + th_i|} - \frac{a|x_i^*|}{1+a|x_i^*|}\right).$$

Therefore, dividing by $t > 0$ on both the sides of inequality (34) and letting $t \to 0$ yield

$$2\langle Ax^* - b, Ah\rangle + \lambda\sum_{i\in\text{supp}(x^*)}\frac{ah_i\text{sgn}\left(x_i^*\right)}{\left(1+a|x_i^*|\right)^2} \geq 0.$$

Obviously, the above-mentioned inequality also holds for $-h$, which leads to equality (11).

2) If $\text{supp}(h) \subseteq \text{Csupp}(x^*)$, then for all $t \in \mathbb{R}$

$$P_a(x^* + th) - P_a(x^*) = \sum_{i\in\text{Csupp}(x^*)}\frac{a|th_i|}{1+a|th_i|}.$$

Hence, it follows from inequality (34) that:

$$t^2\|Ah\|_2^2 + 2t\langle Ax^* - b, Ah\rangle + \lambda\sum_{i\in\text{Csupp}(x^*)}\frac{a|th_i|}{1+a|th_i|} \geq 0.$$

Therefore, for all $t > 0$, we have

$$\begin{aligned}2|\langle Ax^* - b, Ah\rangle| &\leq \|Ah\|_2^2t + \frac{\lambda}{t}\sum_{i\in\text{Csupp}(x^*)}\frac{a|th_i|}{(1+a|th_i|)}\\ &\leq \|Ah\|_2^2t + \frac{\lambda\|h\|_0}{t}.\end{aligned}$$

Thus

$$\begin{aligned}2|\langle Ax^* - b, Ah\rangle| &\leq \min_{t>0}\|Ah\|_2^2t + \frac{\lambda\|h\|_0}{t}\\ &= 2\sqrt{\|Ah\|_2^2\lambda\|h\|_0}\end{aligned}$$

as claimed.

*Proof of Lemma 6*:

1) Let $\text{supp}(h) \subseteq \text{supp}(x^*)$. Then, incorporating equality (11) into inequality (34) yields that, for all $t \in \mathbb{R}$

$$t^2\|Ah\|_2^2 - \lambda\sum_{i\in\text{supp}(x^*)}\frac{tah_i\text{sgn}\left(x_i^*\right)}{\left(1+a|x_i^*|\right)^2} + \lambda(P_a(x^* + th) - P_a(x^*)) \geq 0$$

or equivalently

$$\|Ah\|_2^2 \geq \frac{\lambda}{t^2}\left(\sum_{i\in\text{supp}(x^*)}\frac{tah_i\text{sgn}\left(x_i^*\right)}{\left(1+a|x_i^*|\right)^2} - (P_a(x^* + th) - P_a(x^*))\right). \tag{35}$$

Hence, letting $t \to 0$ on the right-hand side of the above-mentioned inequality, we have inequality (15).

2) If we replace $h$ in inequality (35) with the base vector $e_i$ for every $i \in \text{supp}(x^*)$, then we have the component-wise inequality

$$\|a_i\|_2^2 \geq \frac{\lambda}{t^2}\left(\frac{at\text{sgn}\left(x_i^*\right)}{\left(1+a|x_i^*|\right)^2} - \frac{a|x_i^* + t|}{1+a|x_i^* + t|} + \frac{a|x_i^*|}{1+a|x_i^*|}\right).$$

Particularly, the above-mentioned inequality is available for $t = -x_i^*$. Therefore, we have

$$\|a_i\|_2^2 \geq \frac{\lambda}{x_i^2}\left(-\frac{a|x_i^*|}{(1+a|x_i^*|)^2} + \frac{a|x_i^*|}{1+a|x_i^*|}\right).$$

It follows that:

$$\|a_i\|_2^2 \geq \frac{\lambda a^2}{(1+a|x_i^*|)^2}.$$

From the above-mentioned inequality, inequality (16) immediately follows.

*Proof of Theorem 3*: Assume that $\lambda \in (\|b\|_2^2, \bar{\lambda})$ and (18) holds. If the optimal solution $x^\lambda$ to $(FP_a^\lambda)$ is not the optimal solution to $(FP_a)$, then $Ax^\lambda = b^\lambda \neq b$. Let $y^\lambda$ be the sparsest solution to $Ay = b - b^\lambda$. Then, $A(x^\lambda + y^\lambda) = b$, $\|y^\lambda\|_0 = k \leq m$, and the column submatrix $B^*$ of $A$ consisting of the columns indexed by the set of $\text{supp}(y^\lambda)$ is full column rank. By the definition of $\sigma_{\min}$

$$\sigma_{\min} \leq \frac{\|B^* y^\lambda\|_2^2}{\|y^\lambda\|_2^2}. \tag{36}$$

Using the inequality of the matrix-norm, we obtain

$$\begin{aligned}\|y^\lambda\|_2^2 &\geq \frac{1}{\|A\|_2^2}\|Ay^\lambda\|_2^2 \\ &= \frac{1}{\|A\|_2^2}\|Ax^\lambda - b\|_2^2 \\ &\geq \frac{1}{\|A\|_2^4}\|A^T(Ax^\lambda - b)\|_2^2 \\ &\geq \frac{\lambda^2 a^2}{4\|A\|_2^4}\left(\frac{\lambda - \|b\|_2^2}{\lambda}\right)^4 \end{aligned}\tag{37}$$

where the last inequality holds by inequality (14).

Because $\lambda P_a(y^\lambda) \leq \lambda m$, we have

$$\frac{\lambda P_a(y^\lambda)}{\|y^\lambda\|_2^2} \leq \frac{4m\|A\|_2^4}{\lambda a^2}\left(\frac{\lambda}{\lambda - \|b\|_2^2}\right)^4. \tag{38}$$

Combining it with inequalities (18) and (36), we obtain

$$\frac{\lambda P_a(y^\lambda)}{\|y^\lambda\|_2^2} < \frac{\|B^* y^\lambda\|_2^2}{\|y^\lambda\|_2^2} = \frac{\|Ay^\lambda\|_2^2}{\|y^\lambda\|_2^2}$$

which implies that

$$\lambda P_a(y^\lambda) < \|Ay^\lambda\|_2^2 = \|Ax^\lambda - b\|_2^2.$$

Therefore, we obtain

$$\begin{aligned}&\|A(x^\lambda + y^\lambda) - b\|_2^2 + \lambda P_a(x^\lambda + y^\lambda) \\ &\quad\leq \lambda P_a(x^\lambda) + \lambda P_a(y^\lambda) \\ &\quad< \lambda P_a(x^\lambda) + \|Ax^\lambda - b\|_2^2 \\ &\quad= f(x^\lambda)\end{aligned}\tag{39}$$

which leads to a contradiction that $x^\lambda$ is the optimal solution to $(FP_a^\lambda)$. Hence, the optimal solution to $(FP_a^\lambda)$ also solves $(FP_a)$.

*Proof of Lemma 8*:

1) Define the new variable $\eta = ax + 1$ and substitute it into equality (19); then, the equality can be rewritten as

$$2\eta^3 - 2(1+at)\eta^2 + \lambda a^2 = 0.$$

Due to $t > t_1^*$ and Cardano's root-finding formula expressed in terms of hyperbolic functions (see [32]), the equation has three distinct real roots

$$\eta_0 = \frac{1+at}{3}\left(1 + 2\cos\left(\frac{\phi(t)}{3} - \frac{\pi}{3}\right)\right)$$
$$\eta_1 = \frac{1+at}{3}\left(1 - 2\cos\frac{\phi(t)}{3}\right)$$

and

$$\eta_2 = \frac{1+at}{3}\left(1 + 2\cos\left(\frac{\phi(t)}{3} + \frac{\pi}{3}\right)\right) \tag{40}$$

where

$$\phi(t) = \arccos\left(\frac{27\lambda a^2}{4(1+at)^3} - 1\right).$$

It is obvious that $\eta_0 > \eta_2 > \eta_1$. For $x_i = (\eta_i - 1/a)$, we can also prove $x_0 > x_2 > x_1$. Then, the largest root is $x_0$, i.e., $x_0 = (\eta_0 - 1/a) = g_\lambda(t)$.

2) We set $\eta = 1 - ax$ in equality (20), so $x = 1 - \eta/a$. Then, we can obtain the smallest root with a similar deduction process as the first part

$$x_0 = \frac{1 - \frac{1-at}{3}\left(1 + 2\cos\left(\frac{\phi(t)}{3} - \frac{\pi}{3}\right)\right)}{a}$$

where

$$\phi(t) = \arccos\left(\frac{27\lambda a^2}{4(1-at)^3} - 1\right).$$

Therefore, $x_0 = g_\lambda(t)$ and $x_0 < 0$.

*Proof of Lemma 9*: We discuss $x > 0$, $x = 0$, and $x < 0$ in the following.

1) *$x = 0$:* In this case, $f_\lambda(y) = y^2 + \lambda p_a(|y|)$. It is true that $y^2$ and $\lambda p_a(|y|)$ are increasing with $y > 0$, and they are decreasing with $y < 0$. Thus, $f(0)$ is the least value of $f_\lambda(y)$, i.e., the optimal solution $y^* = 0$ if $x = 0$.

2) *$x > 0$:* It is obvious that $(y - x)^2$ and $\lambda p_a(|y|)$ are decreasing with $y < 0$, so the optimal solution is nonnegative. We just need to consider $y \geq 0$.
In the case $y \geq 0$, we obtain

$$f_\lambda'(y) = 2(y - x) + \frac{\lambda a}{(1+ay)^2}$$

and

$$f_\lambda''(y) = 2 - \frac{2\lambda a^2}{(1+ay)^3}.$$

It is clear that $f_\lambda''(y)$ is increasing. Then, we consider parameter $\lambda$ since it controls the convexity of $f_\lambda(y)$.

a) *$\lambda \leq 1/a^2$:* Because of $\lim_{y\to 0} f_\lambda''(y) = f_\lambda''(0) = 2 - 2\lambda a^2 \geq 0$, $f_\lambda'(y)$ is increasing for $y \geq 0$, and hence, the least value is obtained at $y = 0$ and $f_\lambda'(0) = \lambda a - 2x = 2(\lambda a/2 - x) = 2(t_2^* - x)$.

i) If $0 \leq x \leq t_2^*$, then $f_\lambda(y)$ is positive; thus, the minimum point $y^* = 0$.

ii) If $x > t_2^*$, then $f_\lambda(y)$ is first negative, then positive, and $x > t_2^* > t_1^*$. The minimum point $y^*$ of $f_\lambda(y)$ satisfies

$$f_\lambda'(y^*) = 0 \Leftrightarrow 2y^*(1+ay^*)^2 - 2x(1+ay^*)^2 + \lambda a = 0.$$

Then, the optimal solution is obtained as $y^* = y_0 = g_\lambda(x)$ by Lemma 8.

Specifically, the value of $y^*$ is

$$y^* = \begin{cases} 0, & 0 \leq x \leq t_2^* \\ g_\lambda(x), & x > t_2^* \end{cases}$$

when $\lambda \leq 1/a^2$.

b) $\lambda > 1/a^2$: $f_\lambda'(y)$ is first decreasing and then increasing, and the minimum point of $f_\lambda'(y)$ is $\overline{y} = (\sqrt[3]{\lambda a^2} - 1/a)$. Therefore, the least value

$$f_\lambda'(\overline{y}) = 2(\overline{y} - x) + \frac{\lambda a}{(1+a\overline{y})^2} = 2(t_1^* - x).$$

Then, it is true that $f_\lambda'(y) \geq 2(t_1^* - x)$ with $y \geq 0$.

i) If $0 \leq x \leq t_1^*$, then $f_\lambda(y)$ is increasing with a minimum at $y^* = 0$.

ii) If $x \geq t_2^*$, then $f_\lambda'(0+) \leq 0$, and thus, the function $f_\lambda(y)$ is first decreasing and then increasing, with just one positive optimal point, which is $y^* = g_\lambda(x)$ by Lemma 8.

iii) If $t_1^* < x < t_2^*$, then $f_\lambda'(0+) > 0$, and thus, the function $f_\lambda(y)$ is first increasing, then decreasing, and finally increasing. Hence, $f_\lambda'(y)$ has two positive roots, and the largest root is the minimum point we want. Moreover, the largest root is obtained as $y^* = y_0 = g_\lambda(x)$ by Lemma 8. Thus, we just need to compare $f_\lambda(0)$ with $f_\lambda(y_0)$.

Since a variant of $2(y_0 - x) + (\lambda a/(1+ay_0)^2) = 0$ is $(\lambda a/(1+ay_0)) = 2(x-y_0)(1+ay_0)$, we have

$$\begin{aligned} f_\lambda(y_0) - f_\lambda(0) &= y_0^2 - 2y_0x + \lambda \frac{ay_0}{1+ay_0} \\ &= y_0\left(y_0 - 2x + \frac{\lambda a}{1+ay_0}\right) \\ &= y_0^2(2ax - 1 - 2ay_0) \\ &= 2y_0^2\left(ax - \frac{1}{2} - ay_0\right). \end{aligned}$$

Define a function $\psi(x) = ax - (1/2) - ag_\lambda(x)$.

First, we prove that $x = t_3^*$ is a solution to $\psi(x) = 0$. Due to $\lambda > (1/a^2)$ and $t_3^* = \sqrt{\lambda} - (1/2a) > 0$, there is

$$\cos(\phi(t_3^*)) = \frac{27\lambda a^2}{4(1+at_3^*)^3} - 1 = \frac{27\lambda a^2}{4\left(\frac{1}{2}+a\sqrt{\lambda}\right)^3} - 1.$$

Moreover, we can obtain the following result by formula:

$$\cos(\phi) = 4\cos^3\left(\frac{\phi}{3}\right) - 3\cos\left(\frac{\phi}{3}\right) \quad \left(0 \leq \frac{\phi}{3} \leq \frac{\pi}{3}\right)$$

$$\cos\left(\frac{\phi}{3}\right) = \frac{3\sqrt{8a\sqrt{\lambda}+1} + 4a\sqrt{\lambda} - 1}{4(2a\sqrt{\lambda}+1)}.$$

It is immediate that $g_\lambda(t_3^*) = (\lambda)^{1/2} - 1/a = t_3^* - 1/2a$ after substituting the above-mentioned equation into $g_\lambda(t_3^*)$; therefore, $t_3^*$ is also a solution to $\psi(x)$ in $(t_1^*, t_2^*)$.

Second, we state that function $\psi(x)$ will change its sign at point $x = t_3^*$. We prefer to discuss it in two cases.

*Case 1 [$x \in (t_3^*, t_2^*)$]:* By Lemma 8, we know that $g_\lambda(x)$ is the largest root of the cubic polynomial $2y(1+ay)^2 - 2x(1+ay)^2 + \lambda a = H(y)$ under the condition of $x > t_1^*$. For function $H(y)$, we have $H(x) = \lambda a > 0$ and

$$H\left(x - \frac{1}{2a}\right) = \lambda a - \frac{1}{a}\left(ax + \frac{1}{2}\right)^2.$$

Due to $x > t_3^* = \sqrt{\lambda} - (1/2a)$ and $H(x - (1/2a)) < 0$, there is a root $y = g_\lambda(x)$, such that $g_\lambda(x) \in (x - (1/2a), x)$ for the equation $H(y) = 0$. That is, $x - g_\lambda(x) < (1/2a)$, and thus, $\psi(x) < 0$.

*Case 2 [$x \in (t_1^*, t_3^*)$]:* $H(x - (1/2a)) > 0$ and $H(x) > 0$ hold in this situation. As in Lemma 8, one possible state is that there are two roots $y_0$ and $y_2$ in $(x - (1/2a), x)$. However, we will declare that this is false as follows. With formula (40), there is

$$\begin{aligned} y_0 - y_2 &= \frac{2(1+ax)}{3}\left(\cos\left(\frac{\phi(x) - \pi}{3}\right) - \cos\left(\frac{\phi(x)+\pi}{3}\right)\right) \\ &= \frac{4(1+ax)}{3}\sin\frac{\phi(t)}{3}\sin\frac{\pi}{3}. \end{aligned}$$

Furthermore, $y_0 - y_2 > (1/2a)$ holds for $x > t_1^* > (3/2a^2) - (1/a)$ and $\lambda > (1/a^2)$. This is a contradiction of our assumption that $y_0$ and $y_2$ are in $(x - 1/2a, x)$. Thus, $H(y) = 0$ has no root in $(x - (1/2a), x)$. Therefore, inequality $y_0 = g_\lambda(x) < x - (1/2a)$ holds by $|g_\lambda(x)| \leq |x|$, i.e., $\psi(x) > 0$.

From the discussion earlier, it is true that the optimal solution $y^* = 0$ if $0 < x \leq t_3^*$ and $y^* = y_0 = g_\lambda(x)$ if $x > t_3^*$.

To sum up, we have

$$y^* = \begin{cases} 0, & 0 < x \leq t_3^* \\ g_\lambda(x), & x > t_3^* \end{cases}$$

when $\lambda > (1/a^2)$.

3) $x < 0$: Because

$$\inf_{y\in\mathbb{R}} f_\lambda(y) = \inf_{y\in\mathbb{R}} f_\lambda(-y) = \inf_{y\in\mathbb{R}}\{(y+x)^2 + \lambda p_a(|y|)\}$$

the status of $x > 0$ can be extended to the status of $x < 0$ and formula (21) holds.

According to the results from all cases, the proof is completed.

## ACKNOWLEDGMENT

The authors would like to thank the editor and referees for their comments, which have helped to enrich the content and improve the presentation of the results in this paper.

## REFERENCES

[1] A. M. Bruckstein, D. L. Donoho, and M. Elad, "From sparse solutions of systems of equations to sparse modeling of signals and images," *SIAM Rev.*, vol. 51, no. 1, pp. 34–81, 2009.

[2] E. J. Candés, J. K. Romberg, and T. Tao, "Stable signal recovery from incomplete and inaccurate measurements," *Commun. Pure Appl. Math.*, vol. 59, no. 8, pp. 1207–1223, 2006.

[3] E. J. Candès and T. Tao, "Decoding by linear programming," *IEEE Trans. Inf. Theory*, vol. 51, no. 12, pp. 4203–4215, Dec. 2005.

[4] W. Cao, J. Sun, and Z. Xu, "Fast image deconvolution using closed-form thresholding formulas of $L_q(q = \frac{1}{2}, \frac{2}{3})$ regularization," *J. Vis. Commun. Image Represent.*, vol. 24, no. 1, pp. 31–41, 2013.

[5] R. Chartrand, "Exact reconstruction of sparse signals via nonconvex minimization," *IEEE Signal Process. Lett.*, vol. 14, no. 10, pp. 707–710, Oct. 2007.

[6] S. S. Chen, D. L. Donoho, and M. A. Saunders, "Atomic decomposition by basic pursuit," *SIAM Rev.*, vol. 43, no. 1, pp. 129–159, 2001.

[7] X. Chen, F. Xu, and Y. Ye, "Lower bound theory of nonzero entries in solutions of $\ell_2 - \ell_p$ minimization," *SIAM J. Sci. Comput.*, vol. 32, no. 5, pp. 2832–2852, 2010.

[8] I. Daubechies, M. Defrise, and C. De Mol, "An iterative thresholding algorithm for linear inverse problems with a sparsity constraint," *Commun. Pure Appl. Math.*, vol. 57, no. 11, pp. 1413–1457, Nov. 2004.

[9] I. Daubechies, R. DeVore, M. Fornasier, and C. Güntürk, "Iteratively reweighted least squares minimization for sparse recovery," *Commun. Pure Appl. Math.*, vol. 63, no. 1, pp. 1–38, Jan. 2010.

[10] M. E. Davies and R. Gribonval, "Restricted isometry constants where $l^p$ sparse recovery can fail for $0 < p \leq 1$," *IEEE Trans. Inf. Theory*, vol. 55, no. 5, pp. 2203–2214, May 2009.

[11] D. L. Donoho, "Compressed sensing," *IEEE Trans. Inf. Theory*, vol. 52, no. 4, pp. 1289–1306, Apr. 2006.

[12] D. L. Donoho, "De-noising by soft-thresholding," *IEEE Trans. Inf. Theory*, vol. 41, no. 3, pp. 613–627, May 1995.

[13] D. L. Donoho and X. Huo, "Uncertainty principles and ideal atomic decomposition," *IEEE Trans. Inf. Theory*, vol. 47, no. 7, pp. 2845–2862, Nov. 2001.

[14] D. L. Donoho, "Sparse nonnegative solution of underdetermined linear equations by linear programming," *Proc. Nat. Acad. Sci. USA*, vol. 102, no. 27, pp. 9446–9451, 2005.

[15] D. L. Donoho and M. Elad, "Optimally sparse representation in general (non-orthogonal) dictionaries via $\ell_1$ minimization," in *Proc. Nature Acad. Sci. USA*, vol. 100, Mar. 2003, 2197-2202.

[16] M. Elad, *Sparse and Redundant Representations: From Theory to Applications in Signal and Image Processing*. New York, NY, USA: Springer, 2010.

[17] S. Foucart and M.-J. Lai, "Sparsest solutions of underdetermined linear systems via $\ell_q$-minimization for $0 < q \leq 1$," *Appl. Comput. Harmon. Anal.*, vol. 26, no. 3, pp. 395–407, May 2009.

[18] G. M. Fung and O. L. Mangasarian, "Equivalence of Minimal $\ell_0$ and $\ell_p$ norm solutions of linear equalities, inequalities and linear programs for sufficiently small $p$," *J. Optim. Theory Appl.*, vol. 151, no. 1, pp. 1–10, Oct. 2011. doi: 10.1007/s10957-011-9871-x.

[19] D. Geman and G. Reynolds, "Constrained restoration and the recovery of discontinuities," *IEEE Trans. Pattern Anal. Mach. Intell.*, vol. 14, no. 3, pp. 367–383, Mar. 1992.

[20] T. Goldstein and S. Osher, "The split Bregman method for L1-regularized problems," *SIAM J. Imag. Sci.*, vol. 2, no. 2, pp. 323–343, 2009.

[21] R. Gribonval and M. Nielsen, "Sparse representations in unions of bases," *IEEE Trans. Inf. Theory*, vol. 49, no. 12, pp. 3320–3325, Dec. 2003.

[22] M.-J. Lai and J. Wang, "An unconstrained $l_q$ minimization with $0 < q \leq 1$ for sparse solution of underdetermined linear systems," *SIAM J. Optim.*, vol. 21, no. 1, pp. 82–101, 2011.

[23] H. Li, J. Peng, and S. Yue, "The sparsity of underdetermined linear system via $\ell_p$ minimization for $0 < p < 1$," *Math. Problems Eng.*, vol. 2015, Apr. 2015, Art. no. 584712.

[24] M. Nikolova, "Local strong homogeneity of a regularized estimator," *SIAM J. Appl. Math.*, vol. 61, pp. 633–658, Feb. 2000.

[25] D. Peng, N. Xiu, and J. Yu, "$S_{1/2}$ regularization methods and fixed point algorithms for affine rank minimization problems," *Comput. Optim. Appl.*, vol. 67, no. 3, pp. 543–569, 2017.

[26] J. Peng, S. Yue, and H. Li, "*NP/CMP* equivalence: A phenomenon hidden among sparsity models $\ell_0$ minimization and $\ell_p$ minimization for information processing," *IEEE Trans. Inf. Theory*, vol. 61, no. 7, pp. 4028–4033, Jul. 2015.

[27] R. Saab, R. Chartrand, and O. Yilmaz, "Stable sparse approximations via nonconvex optimization," in *Proc. Int. Conf. Acoust. Speech Signal Process.*, Mar./Apr. 2008, pp. 3885–3888.

[28] Q. Sun, "Recovery of sparsest signals via $\ell_p$ minimization," *Appl. Comput. Harmon. Anaysis*, vol. 32, no. 3, pp. 329–341, 2012.

[29] S. Theodoridis, Y. Kopsinis, and K. Slavakis, "Sparsity-aware learning and compressed sensing: An overview," Nov. 22, 2012. *arXiv:1211.5231*. [Online]. Available: https://arxiv.org/abs/1211.5231

[30] M. Thiao, "Approches de la programmation DC et DCA en data mining," M.S. thesis, INSA-Rouen, Saint-Étienne-du-Rouvray, France, 2011.

[31] M. Wang, W. Xu, and A. Tang, "On the performance of sparse recovery via $l_p$-minimization $(0 \leq p \leq 1)$," *IEEE Trans. Inf. Theory*, vol. 57, no. 11, pp. 7255–7278, Nov. 2011.

[32] F. Xing, "Investigation on solutions of cubic equations with one unknown," *China Acad. J. Electron. Publishing House*, vol. 12, no. 3, pp. 207–218, Jul. 2003.

[33] Z. Xu, X. Chang, F. Xu, and H. Zhang, "$L_{1/2}$ regularization: A thresholding representation theory and a fast solver," *IEEE Trans. Neural Netw. Learn. Syst.*, vol. 23, no. 7, pp. 1013–1027, Jul. 2012.

[34] J. Yang and Y. Zhang, "Alternating direction algorithms for $\ell_1$ problems in compressive sensing," *SIAM J. Sci. Comput.*, vol. 33, no. 1, pp. 250–278, 2011.

[35] W. Yin, S. Osher, D. Goldfarb, and J. Darbon, "Bregman iterative algorithms for $\ell_1$-minimization with applications to compressed sensing," *SIAM J. Imag. Sci.*, vol. 1, no. 1, pp. 143–168, 2008.

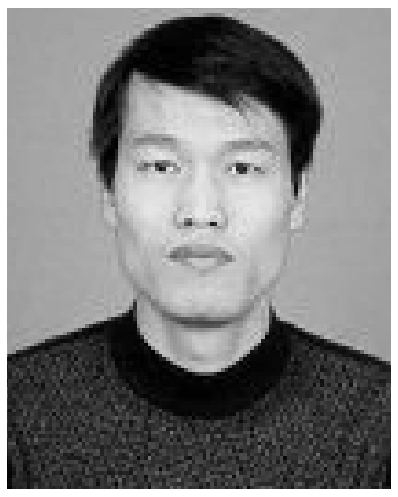

**Haiyang Li** received the B.S. degree in mathematics from Yan'an University, Yan'an, Shaanxi, China, in 1996, and the M.Sc. and Ph.D. degrees in fundamental mathematics from Shaanxi Normal University, Xi'an, China, in 2005 and 2008, respectively.

He is currently a Professor with the School of Mathematics and Information Science, Guangzhou University, Guangzhou, China. His current research interests include quantum logic theory, the machine learning theory, and sparse information processing.

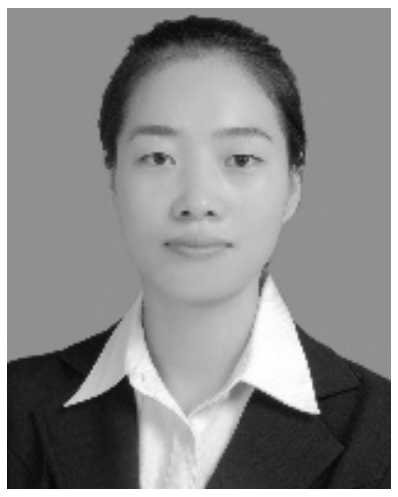

**Qian Zhang** received the B.S. and M.Sc. degrees from the School of Science, Xi'an Polytechnic University, Xi'an, China, in 2015 and 2018, respectively. She is currently pursuing the Ph.D. degree with the School of Computer Science and Engineering, University of Electronic Science and Technology of China, Chengdu, China.

Her current research interests include graph-based semi-supervised learning, Gaussian process, and Bayesian inference.

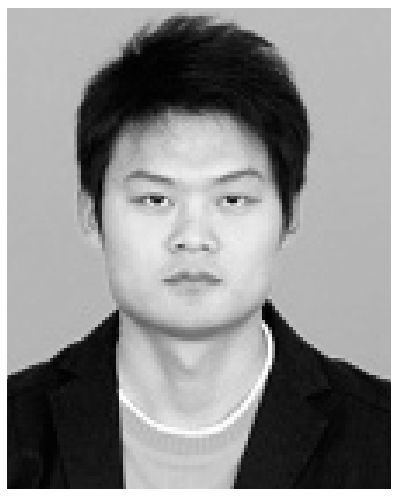

**Angang Cui** received the M.S. degree in applied mathematics from Xi'an Polytechnic University, Xi'an, China, in 2013. He is currently pursuing the Ph.D. degree in applied mathematics and computing mathematics with Xi'an Jiaotong University, Xi'an.

His current research interests include machine learning theory and sparse information processing.

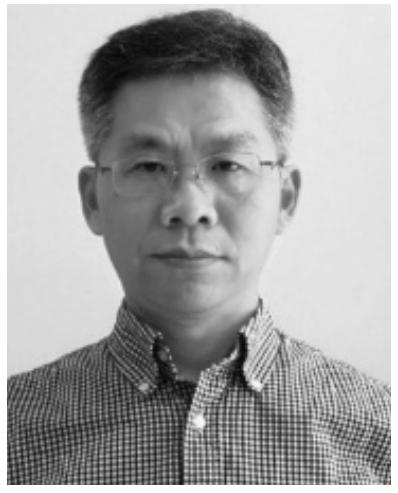

**Jigen Peng** received the B.S. degree in mathematics from Jingxi University, Jingxi, China, in 1989, and the M.Sc. and Ph.D. degrees in applied mathematics and computing mathematics from Xi'an Jiaotong University, Xi'an, China, in 1992 and 1998, respectively.

He is currently a Professor with the School of Mathematics and Information Science, Guangzhou University, Guangzhou, China. His current research interests include nonlinear functional analysis and applications, data set matching theory, the machine learning theory, and sparse information processing.